\def\timestamp{%
Time-stamp: <cosmic-dementia.tex: Tuesday 15-11-2005 at 16:47:31 (cet)>}
\def\stripname Time-stamp: <#1 #2>{#2}
\edef\filedate{\expandafter\stripname\timestamp}

\documentclass%[a4paper]
               {amsart}

\usepackage{amsfonts}
\usepackage{amsrefs}
\usepackage{exscale}

\DeclareMathSymbol\restr\mathbin{AMSa}{"16}
\DeclareMathSymbol\thepo3{AMSa}{"34}
\DeclareMathSymbol\le3{AMSa}{"36}
\DeclareMathSymbol\ge3{AMSa}{"3E}

\DeclareMathSymbol\notsubseteq3{AMSb}{"2A}

\DeclareMathSymbol\I0{AMSb}{`I}
\DeclareMathSymbol\N0{AMSb}{`N}
\DeclareMathSymbol\PP0{AMSb}{`P}
\DeclareMathSymbol\Q0{AMSb}{`Q}
\DeclareMathSymbol\R0{AMSb}{`R}

\newcommand{\cee}{\mathfrak{c}}

\newcommand{\card}[1]{\mathopen|#1\mathclose|}
\newcommand{\norm}[1]{\mathopen\|#1\mathclose\|}
\newcommand{\orpr}[2]{\langle#1,#2\rangle}

\newcommand{\cl}{\operatorname{cl}}
\newcommand{\Fr}{\operatorname{Fr}}
\newcommand{\dom}{\operatorname{dom}}
\newcommand{\ind}{\operatorname{ind}}
\newcommand{\Ind}{\operatorname{Ind}}
\newcommand{\Int}{\operatorname{int}}
\newcommand{\supp}{\operatorname{supp}}
\newcommand{\0}{\mathbf{0}}

\newcommand\QQ{\mathcal{Q}}

\renewcommand\epsilon{\varepsilon}

\newtheorem{proposition}{Proposition}[section]

\newenvironment{numerate}{\enumerate}{\endenumerate}

\begin{document}

\title{Cosmic dimensions}

\author{Alan Dow}
\address{Department of Mathematics\\
         UNC-Charlotte\\
         9201 University City Blvd.\\
         Charlotte, NC 28223-0001}
\email{adow@uncc.edu} 
\urladdr{http://www.math.uncc.edu/\~{}adow}
\thanks{The first author acknowledges support provided by NSF grant
        DMS-0103985.}

\author{Klaas Pieter Hart}
\address{Faculty of Electrical Engineering, Mathematics, and
         Computer Science\\
         TU Delft\\
         Postbus 5031\\
         2600~GA {} Delft\\
         the Netherlands}
\email{K.P.Hart@EWI.TUDelft.NL}
\urladdr{http://aw.twi.tudelft.nl/\~{}hart}

\date{\filedate}

\begin{abstract}
Martin's Axiom for $\sigma$-centered partial orders implies that there
is a cosmic space with non-coinciding dimensions.  
\end{abstract}
\maketitle

\section*{Introduction}

A fundamental result in dimension theory states that the three basic dimension
functions, $\dim$, $\ind$ and $\Ind$, 
coincide on the class of separable metrizable spaces.
Examples abound to show that this does not hold in general outside 
this class.
In \cite{MR0428244} Arkhangel$'$ski\u{\i} asked whether the dimension functions
coincide on the class of cosmic spaces.
These are the regular 
\textbf{co}ntinuous images of \textbf{s}eparable \textbf{m}etric
spaces and they are characterized by the conjunction of regularity and having
a countable network, see~\cite{MR0206907}.
A \emph{network} for a topological space is a collection of (arbitrary)
subsets such that every open set is the union of some subfamily of that
collection.
In~\cite{MR0002526} Vedenisoff proved that $\ind$ and $\Ind$ coincide
on the class of perfectly normal Lindel\"of spaces, 
see also~\cite{MR1363947}*{Section~2.4}.
As the cosmic spaces belong to this class Arkhangel$'$ski\u{\i}'s question
boils down to whether $\dim=\ind$ for cosmic spaces.

In~\cite{MR1661301} Delistathis and Watson constructed, assuming the
Continuum Hypothesis, a cosmic space~$X$ with $\dim X=1$ and $\ind X\ge2$;
this gave a consistent negative answer to Arkhangel$'$ski\u{\i}'s question.

The purpose of this paper is to show that the example can also be constructed
under the assumption of Martin's Axiom for $\sigma$-centered partial orders.
The overall strategy is that of~\cite{MR1661301}: we refine the Euclidean 
topology of a one-dimensional subset~$X$ of the plane to get a topology~$\tau$
with a countable network, such that $\dim(X,\tau)=1$ and in which the boundary
of every non-dense open set is (at least) one-dimensional, so that 
$\ind(X,\tau)\ge2$.
The latter is achieved by ensuring that every such boundary contains a 
topological copy of the unit interval or else a copy of the Cantor set
whose subspace topology is homeomorphic to Kuratowski's graph topology,
as defined in~\cite{Zbl0005.05701}.

The principal difference between our approach and that of~\cite{MR1661301}
lies in the details of the constructions.
In~\cite{MR1661301} the topology is introduced by way of resolutions;
however, some of the arguments given in the paper need emending because,
for example, 
Kuratowski's function does not have the properties asserted and used
in Lemmas~2.2 and~2.3 of~\cite{MR1661301} respectively.
We avoid this and use the Tietze-Urysohn theorem to extend Kuratowski's 
function to the whole plane and thus obtain, per Cantor set, a separable
metric topology on the plane that extends the graph topology.

Also, in~\cite{MR1661301} the construction of the Cantor sets is entwined
with that of the topologies, which lead to some rather inaccessible lemmas.
We separate the two strands and this, combined with the use of partial
orders leads to a cleaner and more perspicuous 
construction of the Cantor sets.
 
We begin, in Section~\ref{sec.thefunction},
with a description of Kuratowski's function. 
We then show how to transplant the graph topology to an arbitrary Cantor
set in the plane.
The remainder of the paper is devoted to a recursive construction of the
necessary Cantor sets and finishes with a verification of the properties of 
the new topology.
An outline of the full construction can be found in Section~\ref{sec.plan}.

\section{Kuratowski's function}
\label{sec.thefunction}

In this section we give a detailed description of Kuratowski's function
(\cite{Zbl0005.05701}, see also~\cite{MR1363947}*{Exercise~1.2.E})
and the resulting topology on the Cantor set.
We do this to make our note self-contained and because the construction
makes explicit use of this description.
We leave the verification of most of the properties to the reader.

Let $C$ be the Cantor set, represented as the topological product~$2^\N$, 
and for $x\in C$ write $\supp x=\{i:x(i)=1\}$.
We let $D$ be the set of~$x$ for which $\supp x$ is finite, partitioned into
the sets $D_k=\{x:\card{\supp x}=k\}$; 
put $k_x=\card{\supp x}$ and $N_x=\max\supp x$ for $x\in D$.
Note that $D_0=\{\0\}$, where $\0$ is the point with all coordinates~$0$.
Let $E=C\setminus D$, the set of~$x$ for which $\supp x$~is infinite.

For $x\in C$ let $c_x$ be the counting function of~$\supp x$, 
so $\dom c_x=\{1,\ldots, k_x\}$ if $x\in D$ and $\dom c_x=\N$ if $x\in E$.
Note that $N_\0=\dom c_\0=\emptyset$.

Now define
$$
f(x)=\sum_{j\in\dom c_x} (-1)^{c_x(j)}2^{-j}
$$
Thus we use the parity of $c_x(j)$ to decide whether to add or 
subtract~$2^{-j}$.
By convention an empty sum has the value~$0$, so $f(\0)=0$.

Notation: if $x\in C$ and $n\in\N$ then $x\restr n$ denotes the restriction
of~$x$ to the set $\{1,2,\ldots,n\}$.
Also, $[x\restr n]$ denotes the $n$-th basic open set around~$x$:
$[x\restr n]=\{y:y\restr n=x\restr n\}$. 

For $x\in D$ we write $V_x=[x\restr N_x]$.
Using the $V_x$ it is readily seen that the sets $D_k$ are relatively
discrete: simply observe that $V_x\cap\bigcup_{i\le k_x}D_i = \{x\}$.
In fact, for a fixed~$k$ the family $\mathcal{D}_k=\{V_x:x\in D_k\}$
is pairwise disjoint.
For later use we put $D_x=\{y\in D_{k_x+1}: y\restr N_x = x\restr N_x\}$
and we observe that $V_x=\{x\}\cup\bigcup\{V_y:y\in D_x\}$.

\subsection{Continuity}

We begin by identifying the points of continuity of $f$.

\begin{proposition}\label{prop.f.cts.at.E}
The function $f$ is continuous at every point of~$E$.\qed
\end{proposition}

The function $f$ is definitely not continuous at the points of~$D$.
This will become clear from the following discussion on the distribution
of the values of~$f$.

\begin{proposition}\label{prop.crowded fibers}
Let $t\in[-1,1]$.
The preimage $f^\gets(t)$ is uncountable, crowded and its intersection
with~$E$ is closed in~$E$. \qed
\end{proposition}

\begin{proposition}\label{prop.accumulate.at.x}
Let $x\in D$ and $k=k_x$.
Then $x$ is an accumulation point of $f^\gets(t)$ if and only if
$f(x)-2^{-k}\le t\le f(x)+2^{-k}$. \qed
\end{proposition}

\subsection{The dimension of the graph}

We identify $f$ with its graph in $C\times[-1,1]$ and
we write $\I=[-1,1]$.
For $x\in D$ we let $I_x=[f(x)-2^{-k_x},f(x)+2^{-k_x}]$.
The discussion in the previous subsection can be summarized by saying that
the closure of~$f$ in~$C\times\I$ is equal to the set
$K= f\cup\bigcup_{x\in D}\bigl(\{x\}\times I_x\bigr)$.

\begin{proposition}
$\ind f\le 1$. \qed  
\end{proposition}

\begin{proposition}\label{prop.ind=0.on.E}
If $x\in E$ then $\ind_{\orpr{x}{f(x)}} f=0$. \qed
\end{proposition}

\begin{proposition}
If $x\in D$ then $\ind_{\orpr{x}{f(x)}} f=1$. \qed
\end{proposition}

We put $\tau_f=\{O_f:O$ open in $C\times\I\}$, where 
$O_f=\{x:\orpr{x}{f(x)}\in O\}$; this is the topology of the graph,
transplanted to~$C$.

\section{Making one Cantor set}
\label{sec:one.cantor.set}

We intend to copy the topology~$\tau_f$ to many Cantor sets in the plane,
or rather, we intend to construct many Cantor sets and copy~$\tau_f$ to
each of them.
Here we describe how we will go about constructing just one Cantor set~$K$,
together with a homeomorphism~$h:C\to K$, and how to refine the topology of 
the plane so that all points but those of~$h[D]$ retain their usual
neighbourhoods and so that at the points of~$h[D]$ the dimension of~$K$ 
will be~$1$.

All we need to make a Cantor set are two maps $\sigma:D\to\R^2$ and 
$\ell:D\to\omega$. 
Using these we define $W(d)=B(\sigma(d),2^{-\ell(d)})$
and $U(d)=B(\sigma(d),2^{-\ell(d)-1})$ for each~$d\in D$.
We want the following conditions fulfilled:
\begin{numerate}
\item the sequence $\langle \sigma(e):e\in D_d\rangle$ converges 
      to~$\sigma(d)$, for all~$d$; \label{cantor.one}
\item $\cl W(e)\subseteq U(d)\setminus \{\sigma(d)\}$ whenever $e\in D_d$;
      \label{cantor.two}
\item $\{\cl W(d):d\in D_n\}$ is pairwise disjoint for all~$n$.
      \label{cantor.three}
\end{numerate}
The following formula then defines a Cantor set:
\begin{equation}\def\theequation{\ddag}
\csname @namedef\endcsname{@currentlabel}{\ddag}
\label{eq:make.Cantor.set}
K=\bigcap_{n=0}^\infty\cl\bigl(\bigcup\{W(d):d\in D_n\}\bigr).
\end{equation}
One readily checks that $\{\sigma(d):d\in D\}$ is a dense subset and that
the map~$\sigma$ extends to a homeomorphism $h:C\to K$ with the property that
$h[V_d]=K\cap W(d)$ for all~$d\in D$.
Also note that in~(\ref{eq:make.Cantor.set}) we could have used the $U(d)$
instead of the $W(d)$ and that even $h[V_d]=K\cap U(d)$ for all~$d$.

Copying the Kuratowski function from~$C$ to~$K$ is an easy matter: 
we let $f_K=f\circ h^{-1}$.
To copy the topology $\tau_f$ to~$K$ \emph{and} to preserve as much as possible
of the Euclidean topology we use the sets~$U(d)$ and~$W(d)$.

We apply the Tietze-Urysohn theorem to extend $f_K$ to a function $\bar f_K$
defined on the whole plane that is continuous everywhere except
at the points of $\sigma[D]$.
The topology~$\tau_K$ that we get by identifying the plane with the 
graph of~$\bar f_K$ is separable and metrizable and its restriction
to $K$ is the graph topology.

As will become clear below we cannot take just any extension of~$f_K$
because we will have to have some amount of continuity at the points
of $\sigma[D]$.
To this end we define for each $d\in D$ a closed  set $F(d)$
by $F(d)=\cl U(d)\setminus\bigcup_{e\in D_d} W(e)$.
The family $\{F(d):d\in D\}$ is pairwise disjoint: if $F(d_1)$ and $F(d_2)$
meet then so do $U(d_1)$ and $U(d_2)$.
Because of conditions (\ref{cantor.two}) and (\ref{cantor.three}) 
above this is only possible if, say, $U(d_1)\supseteq U(d_2)$.
But, unless $d_1=d_2$, this would entail $U(d_2)\subseteq U(e)$ for 
some $e\in D_d$ and so $F(d_2)$ would be disjoint from $F(d_1)$ after all.

The set $K^+=K\cup\bigcup_{d\in D} F(d)$ is closed and we can extend $f_K$
to $K^+$ by setting $f^+_K(x)=f(d)$, whenever $x\in F(d)$.
Because for every $\epsilon>0$ there are only finitely many $d$ for which the
diameter of~$F(d)$ is larger than~$\epsilon$ this extended function
is continuous at all points of~$K\setminus\sigma[D]$.
The new function~$f^+_K$ is certainly continuous at the points 
of~$K^+\setminus K$ (it is even locally constant there), 
so we can apply the Tietze-Urysohn theorem to find a function
$\bar f_K:\R^2\to[-1,1]$ that extends $f^+_K$ and that is continuous
at all points, except those of $\sigma[D]$.

In fact, it not hard to verify that, if $L$ is a subset of the plane
that meets only finitely many of the sets $W(d)$ then the restriction
of $f_K$ to $L$ is continuous.
Indeed, we only have to worry about points in $\sigma[D]$.
But if $d\in D$ then $F(d)\cap L$ contains a neighbour of $\sigma(d)$
in~$L$ and $f^+_K$ is constant on~$F(d)$.

\section{The plan}
\label{sec.plan}

In this section we outline how we will construct a cosmic
topology~$\tau$ on a subset~$X$ of the plane that satisfies 
$\dim(X,\tau)=1$ and $\ind(X,\tau)\ge2$.

We let $\QQ$ denote the family of all non-trivial line segments in the 
plane with rational end points.
Our subset $X$ will be $\R^2\setminus A$, where
$A=\{\orpr{p+\sqrt2}{q}:p,q\in\Q\}$.
Note that $A$~is countable, dense and disjoint from $\bigcup\QQ$.
Also note that, with respect to the Euclidean topology~$\tau_e$, 
one has $\ind(X,\tau_e)=1$:
on the one hand basic rectangles with end~points in~$A$ have 
zero-dimensional boundaries (in~$X$), so that $\ind(X,\tau_e)\le1$, 
and on the other hand, because $X$~is connected we have $\ind(X,\tau_e)\ge1$.

We will construct $\tau$ in such a way that its restrictions 
to $X\setminus\bigcup\QQ$ and each element of~$\QQ$ will be the same 
as the restrictions of~$\tau_e$; 
this ensures that $(X,\tau)$ has a countable network:
take a countable base~$\mathcal{B}$ for the Euclidean topology 
of~$X\setminus\bigcup\QQ$, then
$\QQ\cup\mathcal{B}$ is a network for~$(X,\tau)$.
Also, the $\tau_e$-interior of every open set in $(X,\tau)$ will be nonempty
so that $\bigcup\QQ$ and~$X\setminus\bigcup\QQ$ will be dense
with respect to~$\tau$.

It what follows $\cl$ will be the closure operator with respect to~$\tau$ and
$\cl_e$~will be the Euclidean closure operator.

\subsection*{The topology}

We let $\{(U_\alpha,V_\alpha):\alpha<\cee\}$ numerate all 
pairs of disjoint open sets in the plane whose union is dense and 
for each $\alpha$ we put $S_\alpha=\cl_e U_\alpha\cap \cl_e V_\alpha$.
We shall construct for each $\alpha$ a Cantor set $K_\alpha$
in $X\cap S_\alpha$, unless there is a~$Q_\alpha\in\QQ$ that is contained
in~$S_\alpha$.
The construction of the~$K_\alpha$ will be as described in 
Section~\ref{sec:one.cantor.set}, so that we will be able to extend~$\tau_e$
to a topology~$\tau_\alpha$ whose restriction to~$K_\alpha$ is a copy
of the topology~$\tau_f$.
For notational convenience we let $I$ be the set of~$\alpha$s for which 
we have to construct~$K_\alpha$ and for $\alpha\in\cee\setminus I$ we 
set $\tau_\alpha=\tau_e$.
As an aside we mention that $\cee\setminus I$ is definitely not empty:
if the boundary of $U_\alpha$~is a polygon with rational vertices then
$\alpha\notin I$.

Thus we may (and will) define, for any subset~$J$ of~$\cee$ a 
topology~$\tau_J$: the topology generated by the subbase 
$\bigcup_{\alpha\in J}\tau_\alpha$.
The new topology~$\tau$ will $\tau_\cee$.

There will be certain requirements to be met (the first was mentioned already):
\begin{numerate}
\item The restriction of~$\tau$ to $X\setminus\bigcup\QQ$ and each $Q\in\QQ$
      must be the same as that of the Euclidean topology;\label{cond.one}
\item Different topologies must not interfere: the restriction of~$\tau$
      to~$K_\alpha$ should be the same as that 
      of~$\tau_\alpha$;\label{cond.two}
\item For each $\alpha$, depending on the case that we are in,
      the set $K_\alpha$ or $Q_\alpha$ must be part of the $\tau$-boundary
      of~$U_\alpha$.\label{cond.three} 
\end{numerate}

If these requirements are met then the topology $\tau$ will be as required.
We have already indicated that (\ref{cond.one}) implies that
it has a countable network.

\subsection*{The Inductive dimensions}

To see that $\ind(X,\tau)\ge2$ we take an element~$O$ of~$\tau$ and show that
its boundary is at least one-dimensional.
There will be an~$\alpha$ such that $\cl_eO=\cl_eU_\alpha$: 
there is $O'\in\tau_e$ such that $O\cap \bigcup\QQ= O'\cap\bigcup\QQ$ and 
we can take~$\alpha$ such that $U_\alpha=\Int\cl_e O'$ and 
$V_\alpha=\R^2\setminus\cl_e U$.
In case $\alpha\in I$ the combination of~(\ref{cond.two}) 
and~(\ref{cond.three}) shows that $\ind\Fr O\ge\ind K_\alpha=1$
and in case $\alpha\notin I$ we use (\ref{cond.one}) and~(\ref{cond.three})
to deduce that $\ind\Fr O\ge\ind Q_\alpha=1$.

\subsection*{The covering dimension}

As $\ind(X,\tau)\ge2$ it is immediate that $\dim(X,\tau)\ge1$.
To see that $\dim(X,\tau)\le1$ we consider a finite open cover~$\mathcal{O}$.
Because $(X,\tau)$ is hereditarily Lindel\"of we find that each element 
of~$\mathcal{O}$ is the union of countably many basic open sets.
This in turn implies that there is a countable set~$J$ such that
$\mathcal{O}\subseteq\tau_J$.
The topology~$\tau_J$ is separable and metrizable and it will suffice to show
that $\dim(X,\tau_J)\le1$.

If $J$~is finite then we may apply the countable closed sum theorem:
$O=X\setminus\bigcup_{\alpha\in J}K_\alpha$ is open, hence an $F_\sigma$-set,
say $O=\bigcup_{i=1}^\infty F_i$.
Each $F_i$ is (at most) one-dimensional as is each $K_\alpha$ and
hence so is $X$, as the union of countably many one-dimensional closed
subspaces.

If $J$~is infinite we numerate it as $\{\alpha_n:n\in\N\}$ and set
$J_n=\{\alpha_i:i\le n\}$.
Then $(X,\tau_J)$ is the inverse limit of the sequence 
$\bigl<(X,\tau_{J_n}):n\in\N\bigr>$,
where each bonding map $i_n:(X,\tau_{J_{n+1}})\to (X,\tau_{J_n})$ is the
identity.
By Nagami's theorem 
(\cite{MR0113214}, see also \cite{MR1363947}*{Theorem~1.13.4})
it follows that $\dim(X,\tau_J)\le1$.

\section{The execution}
\label{sec.execution}

The construction will be by recursion on~$\alpha<\cee$.
At stage $\alpha$, if no~$Q_\alpha$ can be found, we take our cue from
Section~\ref{sec:one.cantor.set} and construct 
maps $\sigma_\alpha:D\to S_\alpha$ and $\ell_\alpha:D\to\omega$, in order 
to use the associated balls  
$W_\alpha(d)=B(\sigma_\alpha(d),2^{-\ell_\alpha(d)})$
in formula~(\ref{eq:make.Cantor.set}) to make the Cantor set~$K_\alpha$.
We also get a homeomorphism $h_\alpha:C\to K_\alpha$ as an extension 
of~$d_\alpha$ and use this to copy Kuratowski's function to~$K_\alpha$:
we set $f_\alpha=f\circ h_\alpha^{-1}$.

We use the procedure from the end of Section~\ref{sec:one.cantor.set}
to construct the topology~$\tau_\alpha$.
We let $U_\alpha(d)=B(\sigma_\alpha(d),2^{-\ell_\alpha(d)-1})$,
put $F_\alpha(d)=\cl U_\alpha(d)\setminus\bigcup_{e\in D_d} W_\alpha(e)$
and define $K^+_\alpha$ and $f^+_\alpha$ as above.
We obtain $\tau_\alpha$ as the graph topology from an 
extension~$\bar f_\alpha$ of~$f^+_\alpha$.

\subsection{The partial order}

We construct $\sigma_\alpha$ and $\ell_\alpha$ by an application
of Martin's Axiom to a partial order that we describe in this subsection.
To save on notation we suppress $\alpha$ for the time being.
Thus, $S=S_\alpha$, $\sigma=\sigma_\alpha$, etc.

To begin we observe that $\bigcup\QQ\cap S$ is dense in~$S$:
if $x\in S$ and $\epsilon>0$ then there are points~$a$ and~$b$ with 
rational coordinates in $B(x,\epsilon)$ that belong to~$U$ and~$V$ 
respectively.
The segment $Q=[a,b]$ belongs to~$\QQ$, is contained 
in~$B(X,\epsilon)$ and meets~$S$.
Actually, $Q\cap S$ is nowhere dense in~$Q$ because no subinterval of~$Q$
is contained in~$S$ --- this is where we use the assumption that no element
of~$\QQ$ is contained in~$S$.
There is therefore even a point~$y$ in $Q\cap S$ that belongs 
to $\cl (Q\cap U)\cap\cl(Q\cap V)$: orient $Q$ so that $a$~is its minimum,
then $y=\inf(Q\cap V)$ is as required.
It follows that the set~$S'$ of those $y\in S$ for which there is 
$Q\in\QQ$ such that $y\in\cl (Q\cap U)\cap\cl(Q\cap V)$ is dense 
in~$S$.
We fix a countable dense subset~$T$ of~$S'$.
We also fix a numeration $\{a_n:n\in\N\}$ of~$A$, the complement
of our set~$A$.

The elements $p$ of our partial order $\PP$ have four components:
\begin{numerate}
\item a finite partial function $\sigma_p$ from $D$ to $T$,
\item a finite partial function $\ell_p$ from $D$ to $\omega$,
\item a finite subset $F_p$ of $\alpha\cap I$,
\item a finite subset $\QQ_p$ of $\QQ$.
\end{numerate}
We require that $\dom \sigma_p=\dom \ell_p$ and
we abbreviate this common domain as $\dom p$.
It will be convenient to have $\dom p$ downward closed in~$D$, by which we 
mean that if $e\in\dom p\cap D_d$ then $d\in\dom p$.

The intended interpretation of such a condition is that $\sigma_p$ and $\ell_p$
approximate the maps~$\sigma$ and $\ell$ respectively;
therefore we also write $W_p(d)=B(\sigma_p(d),2^{-\ell_p(d)})$
and $U_p(d)=B(\sigma_p(d),2^{-\ell_p(d)-1})$.
The list of requirements in Section~\ref{sec:one.cantor.set} must be 
translated into conditions that we can impose on~$\sigma_p$ and~$\ell_p$.
\begin{numerate}
\item $\norm{\sigma_p(e)-\sigma_p(d)}<2^{-N_e}$ whenever $d,e\in\dom p$ 
      are such that $e\in D_d$, this will ensure that
      $\langle \sigma(e):e\in D_d\rangle$ will converge to~$\sigma(d)$;
\item $\cl_e W_p(e)\subseteq U_p(d)\setminus\{\sigma_p(d)\}$ 
       whenever $d,e\in\dom p$ 
      are such that $e\in D_d$; and 
\item for every $n$ the family $\{\cl_e W_p(d):d\in D_n\cap\dom p\}$ 
      is pairwise disjoint.
\end{numerate}

The order on $\PP$ will be defined to make $p$~force that for $\beta\in F_p$ 
and $Q\in \QQ_p$ the intersection $\{\sigma(d):d\in D\}\cap(K_\beta\cup Q)$
is contained in the range of~$\sigma_p$, and even that when $d\notin\dom p$ 
the intersection $\cl_eW(d)\cap(K_\beta\cup Q)$~is empty.
We also want $p$~to guarantee that $K\cap\{a_i:i\le\card{\dom p}\}=\emptyset$.

Before we define the order, however, we must introduce an assumption on 
our recursion that makes our density arguments go through with relatively
little effort; unfortunately it involves a bit of notation.

For $x\in\bigcup\QQ$ set $I_x=\{\beta\in I:x\in\sigma_\beta[D]\}$.
For each $\beta\in I_x$ let $d_\beta=\sigma_\beta^{\gets}(x)$ and write
$D_{x,\beta}=D_{d_\beta}$.
If it so happens that $q\in\PP$ and $x=\sigma_q(d)$ for some~$d\in D$ and
if $e\in D_d\setminus\dom q$ then we must be able to choose an extension~$p$
of~$q$ with $e\in\dom p$, without interfering too much with 
the sets $W_\beta(a)$, where $\beta\in I_x$ and $a\in D_{x,\beta}$.
The following assumption enables us to do this 
(and we will be able to propagate it):
\begin{numerate}
\item[$(*)$] If $x\in\bigcup\QQ$ 
   then for every finite subset~$F$ of~$I_x\cap\alpha$ there is 
   an~$\epsilon>0$ such that the family $\mathcal{W}_{F,\epsilon}= 
   \{\cl_eW_\beta(a):\beta\in F, a\in D_{x,\beta}$ and 
   $\sigma(a)\in B(x,\epsilon)\}$  
   is pairwise disjoint.
\end{numerate}
It is an elementary exercise to verify that in such a case the difference
$B(x,\epsilon)\setminus\bigcup\mathcal{W}_{F,\epsilon}$ is connected.
Assumption~$(*)$ will also be useful when we verify some of the properties
of the topology~$\tau$.

We define $p\thepo q$ if 
\begin{numerate}
\item $\sigma_p$ extends $\sigma_q$ and $\ell_p$ extends $\ell_q$,
\item $F_p\supseteq F_q$ and $\QQ_p\supseteq\QQ_q$,
\item if $d\in\dom p\setminus\dom q$ and $i\le\card{\dom q}$ then
      $a_i\notin\cl_e W_p(d)$.
\item if $d\in\dom p\setminus\dom q$ and 
      $J\in \QQ_q\cup\{K_\beta:\beta\in F_p\}$ then $\cl_e W_p(d)$ 
      is disjoint from~$J$.\label{almost.disjoint}
\item if $d\in\dom q$ and $x=\sigma_q(d)$ and if $e\in\dom p\setminus\dom q$
      is such that $e\in D_d$ then $\cl_eW_p(e)$ is disjoint from
      $\cl_eW_\beta(a)$ whenever $\beta\in F_q\cap A_x$ and $a\in D_{x,\beta}$
      \label{propagate.*}
\end{numerate}

It is clear that $p$ and $q$ are compatible whenever $\sigma_p=\sigma_q$ and
$\ell_p=\ell_q$; as there are only countably many possible $\sigma$s and
$\ell$s we find that $\PP$~is a $\sigma$-centered partial order.

\subsection{Dense sets}

In order to apply Martin's Axiom we need, of course, a suitable family
of dense sets.

\subsubsection*{For $\beta<\alpha$ the set $\{p:\beta\in F_p\}$ is dense}
Given $p$ and $\beta$ extend $p$ by adding $\beta$ to~$F_p$.

\subsubsection*{For $Q\in\QQ$ the set $\{p:Q\in\QQ_p\}$ is dense}
Given $p$ and $Q$ extend $p$ by adding $Q$ to~$\QQ_p$.

\subsubsection*{For $n\in\N$ the set $\{p:\card{\dom p}\ge n\}$ is dense}
This follows from the density of the sets below.

\subsubsection*{For $e\in D$ the set $\{p:e\in\dom p\}$ is dense}
Here is where we use assumption~$(*)$.
Since every~$e\in D$ has only finitely many predecessors with respect
to the relation ``$D_d\ni q$'' it will suffice to consider the case
where $q\in\PP$ and $e\in D_d\setminus\dom q$ for some $d\in\dom q$.

We extend $q$ to a condition~$p$ by setting $F_p=F_q$, $\QQ_p=\QQ_q$, 
$\dom p=\{e\}\cup\dom q$ and by defining $d_p(e)$ and $\ell_p(e)$ as follows.
Let $x=\sigma_q(d)$, put $n=k_e$ and consider 
$H=\bigcup\{\cl_e W_q(a):a\in D_{n+1}\cap\dom p\cap D_d\}$.

Fix $\epsilon_1\le 2^{-N_e}$ so that $B(x,2\epsilon_1)$ is 
disjoint from~$H$, this is possible because of condition~(2) in the definition
of the elements of~$\PP$.
Observe that if we choose $\sigma_p(e)$ and $\ell_p(e)$ in such a way
that $\cl_e W_p(e)\subseteq B(x,\epsilon_1/2)$ then $p$~is an element 
of~$\PP$.

Next, using $(*)$, find $\epsilon_2\le\epsilon_1/2$ that works for
the finite set $F_q\cap I_x$. 
The set $W=\{x\}\cup\bigcup\mathcal{W}_{F,\epsilon_2}$ is closed and
does not separate the ball $B(x,\epsilon_2)$, the set $S$ \emph{does}
separate this ball because the latter meets both~$U$ and~$V$.
Therefore we can find a point~$y$ in
$S\cap B(x,\epsilon_2)\setminus W$; we choose $\delta>0$ so small that
$\cl_e B(y,\delta)\subseteq B(x,\epsilon_2)\setminus W$.

The set $S\cap B(y,\delta)$ separates $B(y,\delta)$, hence it is (at least)
one-dimensional,
The union of the $K_\beta$ (for $\beta\in F_q$) together
with the $Q\cap S$ (for $Q\in\QQ_q$) is zero-dimensional because each
individual set is: each $K_\beta$ is a Cantor set and each $Q\cap S$ is
nowhere dense in $Q$ and hence zero-dimensional.
This means that, finally, we can choose $\sigma_p(e)$ in 
$T\cap B(y,\delta)$ but not in this union and then we take $\ell_p(e)$
so large that $\cl_e W_p(e)$ is a subset of $B(y,\delta)$ minus that union.
Also, at this point we ensure that $a_i\notin\cl_e W_p(e)$ 
for $i\le\card{\dom q}$: this is possible because $\sigma_p(e)\notin A$.

We have chosen $W_p(e)$ to meet requirements (3), (4) and (5)
in the definition of~$p\thepo q$.

\subsection{A generic filter}

Let $G$ be a filter on~$\PP$ that meets all of the above dense sets.
Then $\sigma_\alpha=\bigcup\{\sigma_p:p\in G\}$ 
and $\ell_\alpha=\bigcup\{\ell_p:p\in G\}$
are the sought after maps.
We define $W_\alpha$ and $K_\alpha$ as in Section~\ref{sec:one.cantor.set}.

\subsubsection*{Assumption $(*)$ is propagated}
In verifying this we only have to worry about the points 
in $\sigma_\alpha[D]$ of course.

Therefore let $x\in\sigma_\alpha[D]$ and let $F$ be a finite subset 
of~$I_x\cap\alpha$; we have to find an $\epsilon$ for $F'=F\cup\{\alpha\}$.
First fix $\epsilon_1$ that works for $F$ itself.
Next take $p\in G$ such that $d_\alpha\in\dom p$ and $F\subseteq F_p$.
Using condition~(\ref{propagate.*}) in the definition of~$\thepo$ and a density
argument we find that $\cl_e W_\alpha(e)$ is disjoint from
$\cl_e W_\beta(a)$ whenever $e\in D_{x,\alpha}\setminus\dom p$,
$\beta\in F$ and $a\in D_{x,\beta}$.
Now choose $\epsilon$ smaller than $\epsilon_1$ and all
distances $\norm{x-\sigma_\alpha(e)}$, where $e\in D_{x,\alpha}\cap\dom p$.
Then $\mathcal{W}_{F',\epsilon}$ is pairwise disjoint.

\subsubsection*{$K_\alpha$ meets $K_\beta$ in a finite set whenever 
$\beta<\alpha$}

Let $\beta\in\alpha\cap I$ and take $p\in G$ such that $\beta\in F_p$.
Choose $n$ such that $\dom p\subseteq\bigcup_{k\le n}D_k$.
By formula~(\ref{eq:make.Cantor.set}) we know that 
$K_\alpha\subseteq\cl_e\bigl(\bigcup\{W_d:d\in D_{n+1}\}\bigr)$
the latter closure is equal to 
$\bigcup_{k\le n}D_k\cup\bigcup\{W_d:d\in D_{n+1}\}$
and the intersection of \emph{this} set with $K_\beta$ is contained 
in $\dom p$; this follows from condition~\ref{almost.disjoint} in the
definition of~$\thepo$.

\subsubsection*{$K_\alpha$ meets each $Q\in\QQ$ in a finite set}

The proof is identical to the previous one: take $p\in G$ with
$Q\in\QQ$.

\section{The remaining properties of the topologies}

We check conditions (\ref{cond.one}), (\ref{cond.two}) 
and (\ref{cond.three}) from Section~\ref{sec.plan}.

A useful observation is that a typical new basic neighbourhood of a point~$x$
of~$\bigcup\QQ$ contains a set of the form 
$O(x,\epsilon,G)=B(x,\epsilon)\cap\bigcap_{\beta\in G}F_\beta(d_\beta)$,
where $\epsilon>0$ and $G$~is a finite subset of~$I_x$.

\subsection{$X\setminus\bigcup\QQ$ retains its Euclidean topology}

This is immediate from the observation that every function $\bar f_\alpha$
(for $\alpha\in I$) is continuous at the points of~$X\setminus\bigcup\QQ$.

\subsection{Each $Q\in\QQ$ retains its Euclidean topology}

We should show that $\bar f_\alpha\restr Q$ is continuous for each~$\alpha$
in~$I$ and each $Q\in\QQ$.
The only points at which this restriction could possibly be discontinuous
are those in $\sigma_\alpha[D]\cap Q$, which is a finite set.
Let $d\in D$ be such that $x=\sigma_\alpha(d)\in Q$.
By construction all but finitely many of the sets $\cl W_\alpha(e)$,
where $e\in D_d$, meet $Q$.
This implies that $F_\alpha(d)\cap Q$ is actually a neighbourhood
of $x$ in~$Q$.
As $\bar f_\alpha$ is constant on~$F_\alpha(d)$ this shows that 
$\bar f_\alpha\restr Q$ is continuous at~$x$.

\subsection{$\tau_\alpha$ and $\tau_\beta$ do not interfere}

If $\alpha\neq\beta$ then there are only finitely many points 
in~$K_\alpha\cap K_\beta$ and it is only at these points that
$\tau_\alpha$ and $\tau_\beta$ might interfere and even then only
at a point of $\sigma_\alpha[D]\cap\sigma_\beta[D]$.
Let $x$ be such a point and apply assumption~$(*)$ to the 
set $G=\{\alpha,\beta\}$ to find $\epsilon>0$ such that 
$\mathcal{W}_{G,\epsilon}$ is pairwise disjoint.
But then $f^+_\alpha\restr K_\beta$ is constant on a neighbourhood of~$x$
in~$K_\beta$, namely $O(x,\epsilon,\{\alpha\})\cap K_\beta$ and, by symmetry,
$f^+_\beta\restr K_\alpha$ is constant on the neighbourhood 
$O(x,\epsilon,\{\beta\})\cap K_\alpha$ of~$x$ in~$K_\alpha$.

\subsection{$Q_\alpha$ is still in the boundary of $U_\alpha$}

If $x\in Q_\alpha$ then, by construction, all points of the intersection 
$Q_\alpha\cap O(x,\epsilon,G)$ (except $x$ itself) belong to the Euclidean 
interior of $O(x,\epsilon,G)$,
Because these points are in the boundary of~$U_\alpha$ that interior meets
both~$U_\alpha$ and~$V_\alpha$.
Therefore each basic neighbourhood of~$x$ meets these sets as well.

\subsection{$K_\alpha$ is still in the boundary of $U_\alpha$}

Let $x\in K_\alpha$, assume $I_x\neq\emptyset$ and consider some 
$O(x,\epsilon,G)$.

If $\alpha\notin I_x$ then the same argument as above will work:
the intersections $B(x,\epsilon)\cap K_\alpha$ and 
$O(x,\epsilon,G)\cap K_\alpha$ are equal when $\epsilon$ is small enough.

If $\alpha\in I_x$ then we assume $\alpha\in G$ and observe that if $\epsilon$
is small enough then $O(x,\epsilon,G)$ is a Euclidean neighbourhood of 
many points of $S_\alpha$.

\end{document}